\documentclass[12pt]{article}
\usepackage{amssymb,amsmath,amsthm}
\begin{document}
\author{Federico Montecalvo
\\ School of Mathematical Sciences
\\ Queen Mary, University of London
\\ Mile End Road, London E1 4NS U.K.
\\
\\ \texttt{F.Montecalvo@qmul.ac.uk}
}
\title{Constructions of General Covering Designs\footnote{The submitted version of this paper is called \textit{Some Constructions of General Covering Designs} and does not contain Section \ref{sec5_4_6}. Instead, it contains a section where we presented  constructions for $(v,6,5,7)$ covers.}}
\date{December 2012}
\maketitle
%%% ----------------------------------------------------------------------
%\newcommand{\proofend}{\hspace*{\fill}\mbox{$\square$}}
\newtheorem{teok57}{Theorem}
\theoremstyle{definition}\newtheorem{def_block_des}[teok57]{Definition}\theoremstyle{plain}
\theoremstyle{definition}\newtheorem{t_design}[teok57]{Definition}\theoremstyle{plain}
\newtheorem{STS}[teok57]{Theorem}
\theoremstyle{definition}\newtheorem{i_resolvable}[teok57]{Definition}\theoremstyle{plain}
\newtheorem{baker1976}[teok57]{Theorem}
\newtheorem{teirlinck1994}[teok57]{Theorem}
\theoremstyle{definition}\newtheorem{complete_graph}[teok57]{Definition}\theoremstyle{plain}
\theoremstyle{definition}\newtheorem{res_cov}[teok57]{Definition}\theoremstyle{plain}
\newtheorem{res_cov_q_r}[teok57]{Theorem}
\newtheorem{res_cov_q_r_upp}[teok57]{Theorem}
\newtheorem{gen_splicing}[teok57]{Theorem}
\newtheorem{spec_splicing}[teok57]{Theorem}
\newtheorem{new34-6}[teok57]{Theorem}
\newtheorem{teo63m}[teok57]{Theorem}
\newtheorem{genteo63m}[teok57]{Theorem}
\newtheorem{genteo634}[teok57]{Corollary}
\newtheorem{teok46}[teok57]{Theorem}
\theoremstyle{remark}\newtheorem{remarkupper646}[teok57]{Remark}\theoremstyle{plain}
\newtheorem{teok45}[teok57]{Theorem}
\newtheorem{teoupperformulak46}[teok57]{Theorem}
\newtheorem{teo546}[teok57]{Theorem}
%\newtheorem{teo657}[teok57]{Theorem}
%\newtheorem{teo29_657}[teok57]{Theorem}
%%% ----------------------------------------------------------------------
\begin{abstract}
Given five positive integers $v, m,k,\lambda$ and $t$ where $v \geq k \geq t$
and $v \geq m \geq t,$ a $t$-$(v,k,m,\lambda)$ \textit{general
covering design} is a pair $(X,\mathcal{B})$ where $X$ is a set
of $v$ elements (called points) and $\mathcal{B}$ a multiset of
$k$-subsets of $X$ (called blocks) such that every $m$-subset of $X$
intersects (is covered by) at least $\lambda$ members of $\mathcal{B}$ in at least
$t$ points.

In this article we present new constructions for general covering designs and we generalize some others. By means of these constructions we will be able to obtain some new upper bounds on the minimum size of such designs.

\end{abstract}
  % ------------------------------------------------------------------------
\section{Introduction}
Given five positive integers $v, m,k,\lambda$ and $t$ where $v \geq k \geq t$
and $v \geq m \geq t,$ a $t$-$(v,k,m,\lambda)$ \textit{general covering design} (or \textit{general cover}) is a pair $(X,\mathcal{B})$ where $X$ is a set
of $v$ elements (called points) and $\mathcal{B}$ a multiset of
$k$-subsets of $X$ (called blocks) such that every $m$-subset of $X$
intersects (is covered by) at least $\lambda$ members of $\mathcal{B}$ in at least
$t$ points.

It is easy to verify that a $t$-$(v,k,m,\lambda)$ general cover is also a $(t-1)$-$(v,k,m-1,\lambda)$ general cover. A $t$-$(v,k,m,\lambda)$ general covering design $(X,\mathcal{B})$ is said to be \textit{optimal} if:
\begin{center}
$\vert \mathcal{B}\vert = \textrm{min} \{\vert\mathcal{A}\vert :$ there is a $t$-$(v,k,m,\lambda)$ general covering design $(X,\mathcal{A})\}$.
\end{center}
In this case, the cardinality of $\mathcal{B}$ is called the \textit{general covering number} and denoted by $C_{\lambda}(v,k,t,m)$.

Given a $t$-$(v,k,m,1)$ general covering design $(X,\mathcal{B})$, the set $\mathcal{C} = \{X \setminus B: B \in \mathcal{B} \}$ is said to be the collection of the coblocks of $(X,\mathcal{B})$ and the pair $(X, \mathcal{C})$ is called the complement of $(X,\mathcal{B})$.

Applications to error-trapping decoding, data compression and lottery systems have led many special cases of general covering designs to be investigated. Let us describe the most studied in the literature:
\begin{list}{}{}
\item  \textbf{Covering Designs}: When $m=t$ and $\lambda = 1$, a $t$-$(v,k,m,\lambda)$ general covering design is said to be a $(v,k,t)$ covering design. The general covering number is simply called \textit{covering number} and denoted by $C(v,k,t)$. There is an extensive literature on covering designs. For an excellent survey please refer to \cite{gordon1,mil,milmul}. Covering designs are applied to error-trapping decoding \cite{chan}. Here the number of the blocks determines the complexity of the decoding procedure. So, optimal covering designs are of special interest.
\item  \textbf{Tur\'an Systems}: When $k=t$ and $\lambda = 1$, a $t$-$(v,k,m,\lambda)$ general covering design is said to be a $(v,k,m)$ Tur\'an system.
The general covering number is called \textit{Tur\'an number} and denoted by $T(v,k,m)$. By taking the coblocks of a $(v,k,t)$ covering design, we always obtain a $(v,v-k,v-t)$ Tur\'an system. Conversely, if we take the coblocks of a $(v,k,m)$ Tur\'an system we always obtain a $(v,v-k, v-m)$ covering design.  Therefore: \mbox{$T(v,k,m) = C(v,v-k,v-m)$} and \mbox{$C(v,k,t) = T(v,v-k,v-t)$}. For a survey please refer to \cite{brower, decaen1, sid}.
\item  \textbf{Lotto Designs}: When $\lambda = 1$, a $t$-$(v,k,m,\lambda)$ general covering design is said to be a $(v,k,t,m)$ lotto design (or \textit{cover}). We will generally use the latter definition in the following sections. The general covering number is called \textit{lotto} (or \textit{cover}) \textit{number} and denoted either by $L(v,k,t,m)$ or by $C(v,k,t,m)$. From the definition, both covering designs and Tur\'an systems can be seen as special cases of lotto designs where $m = t$ and $k=t$ respectively. Therefore $C(v,k,t) = C(v,k,t,t) = T(v,v-k,v-t)$ and $T(v,k,m) = C(v,k,k,m) = C(v,v-k,v-m)$. As the name suggests, lotto designs find application to national lotteries \cite{bert, cre, fur}, but they are also applied to data compression algorithms, as described in \cite{etz1}. Several studies have focused on establishing upper and lower bounds on $C(v,k,t,m)$. Currently, the situation is as
follows:
\begin{list}{-}{}
\item Only for few values of $v,k,t$ and $m$ the cover number $C(v,k,t,m)$
has been found (see \cite{bate,bou,brower,li1}).
\item Constructions and lotto tables have been published in international journals (see \cite{batev, bert, cre, li1}).
\item Upper bounds on $C(v,k,t,m)$ are available on web sites (see \cite{belic, fromme, gordon}).
\item Results on lower bounds have also been published (see \cite{fur, li2}).
\end{list}
\end{list}

General covers should not be confused with a class of objects called \textit{generalized covering designs} which were recently introduced by Bailey et al. in \cite{bai}. Generalized covering designs simultaneously generalize covering designs and covering arrays. For further information and details on this class of objects, the reader is referred to the aforementioned reference.
\section{Background}
In this section we present definitions and known results on design theory which will be used throughout this article.
\begin{def_block_des}
A \textit{block design} is a pair $(X,\mathcal{B})$ such that:
\begin{enumerate}
\item $X$ is a set of elements called \textit{points}.
\item $\mathcal{B}$ is a multiset (\emph{collection}) of non-empty subsets of $X$ called \textit{blocks}.
\end{enumerate}
\end{def_block_des}
The cardinality of $X$ is said to be the \textit{order} of a block design $(X,\mathcal{B})$. Two block designs $(X,\mathcal{A})$ and $(X,\mathcal{B})$  are called  \textit{disjoint} if $\mathcal{A} \cap \mathcal{B} = \emptyset$.  The \textit{product} of two block designs $(X_1, \mathcal{A})$ and $(X_2, \mathcal{B})$ is defined as $(X_1 \cup X_2,\mathcal{A}\mathcal{B})$ where $\mathcal{A}\mathcal{B} = \{A \cup B: A \in \mathcal{A}, B \in \mathcal{B} \}$.
\begin{t_design}
A \emph{$t$-$(v,k,\lambda)$-design} is a pair $(X, \mathcal{B})$ where $X$ is a set
of $v$ elements (called points) and $\mathcal{B}$ a multiset of
$k$-subsets of $X$ (called blocks) such that every $t$-subset of $X$
is contained in exactly $\lambda$ blocks.
\end{t_design}

The general term \emph{$t$-design} is often used to indicate any $t$-($v,k,\lambda)$-design. When $\lambda = 1$, a $t$-($v,k,1)$-design is often called a \textit{Steiner system} and denoted by $S(v,k,t)$. If $t=2$ and $k=3$, a Steiner system is called a \textit{Steiner triple system} and denoted by STS$(v)$ and if $t=3$ and $k=4$ it is called a \textit{Steiner quadruple system} and denoted by SQS$(v)$.

When $\lambda > 1$, the union of two collections of blocks $\mathcal{A}$ and $\mathcal{B}$ of $t$-designs (or general covering designs) is a multiset union. Therefore, if a block $C$ appears $r_1$ times in $\mathcal{A}$ and $r_2$ times in $\mathcal{B}$, $C$ will appear $\max\{r_1, r_2\}$ times in  $\mathcal{A} \cup \mathcal{B}$.

Many results on the necessary and sufficient conditions for the existence of $t$-designs have been found. Here we report the one on Steiner triple systems:
\begin{STS} {\rm \cite{kirk}}
There exists an $\emph{STS}(v)$ if and only if $v \equiv 1,3 \;(\emph{mod}\;6)$, $v \geq 7$.
\end{STS}
A $t$-design $(X, \mathcal{B})$ is said to be \textit{$\alpha$-resolvable} if there exists a partition of the collection $\mathcal{B}$ into parts called \textit{$\alpha$-parallel classes} (or \textit{$\alpha$-resolution classes}) such that each point of $X$ occurs exactly in $\alpha$ blocks in each class. When $\alpha = 1$, $\alpha$ is omitted. 

Another interesting concept is the one of $i$-partitionable designs:
\begin{i_resolvable}
A Steiner system $S(v,k,t)$ is called \emph{$i$-partitionable}, $0 < i < t$, if the collection of its blocks can be partitioned into Steiner systems $S(v,k,i)$.
\end{i_resolvable}
With regard to $i$-partitionable designs, the following two important theorems hold:
\begin{baker1976}\label{baker1976} {\rm \cite{bak}}
For any positive integer $n$ there exists a 2-partitionable \emph{SQS($4^n$)}.
\end{baker1976}
\begin{teirlinck1994}\label{teirlinck1994}  {\rm \cite{teirl}}
For any positive integer $n$ there exists a 2-partitionable \emph{SQS($2p^n + 2$)}, $p \in \{7,31,127\}$.
\end{teirlinck1994}
When $k = 2$, we often talk in terms of graphs rather than designs.
\begin{complete_graph}
The \emph{complete graph of order $n$}, denoted by $K_n$, is a regular graph with $n$ vertices such that each pair of vertices is an edge.
\end{complete_graph}
The number of edges of the complete graph $K_n$ is $\frac{n(n-1)}{2}$, that is, all the possible pairs of vertices.

A \emph{$1$-factor} of a graph $G$ is a set $E$ of edges such that every vertex of $G$ is incident to exactly one edge of $E$. A \emph{$1$-factorization} is a partition of the edges of a graph  into $1$-factors. In term of designs, a $1$-factorization of the complete graph $K_n$ corresponds to a partition of the Steiner system $S(n,2,2)$ (i.e. the set of all the pairs from $n$) into parallel classes. Clearly, $n$ must be even.

A definition of resolvability can be extended to covering designs as follows:
\begin{res_cov}
A $(v,k,t)$ covering design $(X, \mathcal{B})$ is \emph{resolvable} if $\mathcal{B}$ can be partitioned into parts called \emph{parallel classes} (or \emph{resolution classes}) each of which in turn partitions $X$.
\end{res_cov}
The number of blocks in a parallel class is necessarily $v/k$.

Let $r(q,k)$ denote the minimum number of parallel classes in a resolvable $(kq,k,2)$ covering design. When $q = 1$, $r(q,k)$ is trivially equal to 1. The following results hold:

\begin{res_cov_q_r_upp}\label{haemer_gen} {\rm \cite{vhaemers}} When $q > 1$, $r(q,k) \geq q + 1$. Equality holds if and only if $q$ divides $k$ and $q$ is the order of an affine plane.
\end{res_cov_q_r_upp}
For small values of $q$:
\begin{res_cov_q_r}\label{haemer} {\rm \cite{vhaemers}}
$\;$
\begin{enumerate}
\item $r(2,k) = 3$ if $k$ is even, $4$ if $k$ is odd;
\item $r(3,k) = 4$ if $k \equiv 0 \;(\emph{mod}\;3)$, $5$ otherwise;
\item $r(4,k) = 5$ if $k \equiv 0 \;(\emph{mod}\;4)$, $7$ if $k \in \{2,3\}$, $6$ otherwise.
\end{enumerate}
\end{res_cov_q_r}
Another interesting concept is the one of \textit{large set of coverings}. Given a set $X$ of size $v$ and a positive integer $k$, let $\binom{X}{k}$ be the set of all $k$-subsets of $X$ and let $\mu(v,k)$ denote the minimum number of optimal $(v,k,k-1)$ covering designs $(X, \mathcal{B}_1),(X, \mathcal{B}_2), \ldots, (X, \mathcal{B}_{\mu(v,k)})$ such that $\bigcup_{i=1}^{\mu(v,k)} \mathcal{B}_i = \binom{X}{k}$. Let $\lambda(v,k)$ denote instead the maximum number of disjoint optimal $(v,k,k-1)$ covering designs defined on $X$. Then a large set of coverings is obtained when $\lambda(v,k) = \mu(v,k)$. %The following results for $\mu(10,4)$ hold:
%\begin{large_set_cov} \label{etzio} {\rm \cite{etz}} $\displaystyle \mu(10,4) \leq 10$.
%\end{large_set_cov}

\bigskip
In the following sections, given a partition $X_1, \ldots, X_n$ of a set $X$ of size $v$, a positive integer $m \leq v$, and $n$ positive integers $a_1 \leq \vert X_1 \vert, \ldots, a_n \leq \vert X_n \vert$ such that $\sum_{i=1}^{n} a_i = m$, we will assume that $[a_1,\ldots, a_n]$ denotes the subset of $\binom{X}{m}$ whose elements $M$ satisfy $\vert M \cap X_i \vert = a_i$, for $1 \leq i \leq n$.
\section{Point Splicing Constructions}
Etzion et al. \cite{etz1} described a construction for constant weight covering codes called one-bit splicing. It was actually a construction for $(v,k,m)$ Tur\'an systems. The objective was to start from a Tur\'an system of order $v$ to obtain a Tur\'an system of order $v+1$. In the next section we present a simple generalization: We start from a general covering design of order $v$ to obtain a general covering design of order $v+n$.
\subsection {Point Splicing Construction for $t$-$(v,k,m,\lambda)$ General Covers}\label{n-point}
Let $(X,\mathcal{B})$ be a $t$-$(v,k,m,\lambda)$ general covering design and $n$ be the size of a set $S$ such that $X \cap S = \emptyset$ and $n \leq k-t+1$, where $t > 2$. For every $x \in X$, define $\mathcal{B}(x) = \{B \setminus \{x\}: B \in \mathcal{B}, x \in B\}$. Choose $a \in X$ such that for any $x \in X$ we have $\vert \mathcal{B}(a) \vert \leq \vert \mathcal{B}(x) \vert$. Let $(X \setminus \{a\},\mathcal{C)}$ be a $(t-2)$-$(v-1,k-n-1,m-2,\lambda)$ general covering design and $\mathcal{B}_1$, $\mathcal{B}_2$, and $\mathcal{B}_3$ be three collections of blocks as defined below:
\begin{list}{}{}
\item $\mathcal{B}_1 = \mathcal{B}$.
\item $\mathcal{B}_2 = \{B \cup \{s\}: B \in \mathcal{B}(a), s \in S\}$.
\item $\mathcal{B}_3 = \{C \cup S \cup \{a\} : C \in \mathcal{C}\}$. \end{list}
Our objective is to obtain a $t$-$(v + n,k,m,\lambda)$ general covering design on the set  $(X \cup S)$ and we claim that $(X \cup S, \mathcal{B}_1 \cup \mathcal{B}_2 \cup \mathcal{B}_3)$ meets the objective.
\begin{gen_splicing}\label{gen_splicing}
$(X \cup S, \mathcal{B}_1 \cup \mathcal{B}_2 \cup \mathcal{B}_3)$ is a $t$-$(v + n,k,m,\lambda)$ general covering design.
\end{gen_splicing}
\begin{proof} Let $M$ be an $m$-subset of $X \cup S$:
\begin{list}{}{}
\item If $M \cap (S \cup \{a\}) = \emptyset$ or $\{a\}$, then there exist at least $\lambda$ blocks in $\mathcal{B}_1$ that cover $M$ in $t$ points.
\item If $M \cap (S \cup \{a\}) = \{s\}$, $s \in S$, then there exist at least $\lambda$ blocks in $\mathcal{B}_1 \cup \mathcal{B}_2$ that cover $M$ in $t$ points.
\item If  $\vert M \cap (S \cup \{a\}) \vert = \delta$,
$2 \leq \delta < t$,   then there exist at least $\lambda$ blocks in $\mathcal{B}_3$ that cover $M$ in $t$ points because
$(X \setminus \{a\}, \mathcal{C})$ is a $(t-2)$-$(v-1,k-n-1,m-2,\lambda)$ general covering design and therefore a $(t-\delta)$-$(v-1,k-n-1, m - \delta,\lambda)$ general covering design as well.
\item If  $\vert M \cap (S \cup \{a\}) \vert = \delta \geq t$, then $M$ is clearly covered by each block $B \in \mathcal{B}_3$.
\end{list}
Therefore $(X \cup S, \mathcal{B}_1 \cup \mathcal{B}_2 \cup \mathcal{B}_3)$ is a $t$-$(v + n,k,m,\lambda)$ general covering design. \end{proof}
By counting arguments, $\min_{x \in X} \vert \mathcal{B}(x) \vert \leq \left\lfloor \frac{k\vert \mathcal{B} \vert}{v}\right \rfloor$,  therefore, as a consequence of the construction above:
\begin{displaymath}\label{eqobstwice}
C_{\lambda}(v+n,k,t,m) \leq n \left\lfloor \frac{k}{v}
C_{\lambda}(v,k,t,m)\right \rfloor + C_{\lambda}(v,k,t,m) + C_{\lambda}(v-1,k-n-1,t-2, m-2).
\end{displaymath}
\subsection {Point Splicing Construction for $(v,k,4,6)$ Covers}\label{splice_constr_spec}
We introduce a point splicing construction specific for ($v,k,4,6)$ covers. Similar in spirit to a construction for $(v,4,6)$ Tur\'an systems presented by Etzion et al. \cite{etz1}, it permits us to obtain a $(v+3,k,4,6)$ cover from a $(v,k,4,6)$ cover. Moreover, the technique on which this construction is based allows us to derive a $(v+3,k,4,5)$ cover from a $(v,k,4,5)$ cover.

Let $k \geq 5$. Let $(X,\mathcal{B})$ be a $(v,k,4,6)$ cover.  For every $x \in X$, let $\mathcal{B}(x)$ be defined as in Section \ref{n-point}. Choose $a \in X$ such that for any $x \in X$ we have $\vert \mathcal{B}(a) \vert \leq \vert \mathcal{B}(x) \vert$. Let $b, c$ and $d$ be three new points such that $X \cap \{b,c, d\} = \emptyset$. Let $X_{1,1}, X_{1,2}$, $X_{2,1}, X_{2,2}$, $X_{3,1}, X_{3,2}$ be a partition of $X \setminus \{a\}$.  Let us call a covering design with $t=2$ and block size $k$ a $(2,k)$-\emph{covering}. %a $(v,3,2)$ covering design and a pair-by-quadruple covering denotes a $(v,4,2)$ covering: 
Then take the following covering designs:
\begin{equation*}
\text{Three $(2, k-3)$-coverings}
\begin{cases}
(X_{1,1}\cup X_{2,1}, \mathcal{C}_1)\\
(X_{1,2}\cup X_{2,2}, \mathcal{C}_2)\\
(X_{3,1} \cup X_{3,2}, \mathcal{C}_3)\end{cases}
\end{equation*}
\begin{equation*}
\text{Three $(2, k-3)$-coverings}
\begin{cases}
(X_{1,1} \cup X_{2,2}, \mathcal{D}_1)\\
(X_{1,2} \cup X_{2,1}, \mathcal{D}_2)\\
(X_{3,1} \cup X_{3,2}, \mathcal{D}_3)
\end{cases}
\end{equation*}
\begin{equation*}
\text{Three $(2, k-2)$-coverings}
\begin{cases}
(X_{1,1}\cup X_{1,2}, \mathcal{E}_1)\\
(X_{2,1} \cup X_{2,2}, \mathcal{E}_2)\\
(X_{3,1} \cup X_{3,2}, \mathcal{E}_3)
\end{cases}
\end{equation*}
\newline\indent
%\begin{list}{}{}
%\item Three $(2, k-3)$-coverings $(X_{1,1}\cup X_{2,1}, \mathcal{C}_1)$, $(X_{1,2}\cup X_{2,2}, \mathcal{C}_2)$ and  $(X_{3,1} \cup X_{3,2}, \mathcal{C}_3)$.
%\item Three $(2, k-3)$-coverings $(X_{1,1} \cup X_{2,2}, \mathcal{D}_1)$, $(X_{1,2} \cup X_{2,1}, \mathcal{D}_2)$ and $(X_{3,1} \cup X_{3,2}, \mathcal{D}_3)$.
%\item Three $(2, k-2)$-coverings $(X_{1,1}\cup X_{1,2}, \mathcal{E}_1)$, $(X_{2,1} \cup %X_{2,2}, \mathcal{E}_2)$ and $(X_{3,1} \cup X_{3,2}, \mathcal{E}_3)$.
%\end {list}
The designs above have the following properties:
\begin{enumerate}
\item $(X \setminus \{a\}, \bigcup_{i=1}^{3} \mathcal{C}_i)$ is a $(v-1,k-3,2,4)$ cover. \label{pro1}
\item $(X \setminus \{a\}, \bigcup_{i=1}^{3} \mathcal{D}_i)$ is a $(v-1,k-3,2,4)$ cover. \label{pro2}
\item $(X \setminus \{a\}, \bigcup_{i=1}^{3}\mathcal{E}_i) $ is a $(v-1,k-2,2,4)$ cover. \label{pro3}
\item For any triple $\{x,y,z \} \subset X \setminus \{a\}$, there exists $B \in \bigcup_{i=1}^{3} \mathcal{C}_i \cup \mathcal{D}_i \cup \mathcal{E}_i$  such that $ \vert B \cap  \{x,y, z \} \vert \geq 2$. \label{pro4}
\end{enumerate}
We can now proceed to build a $(v+3,k,4,6)$ cover. Define:
\begin{list}{}{}
\item $\mathcal{B}_1 = \mathcal{B}$.
\item $\mathcal{B}_2 = \{B \cup \{p\}: B \in \mathcal{B}(a), p \in \{b,c,d\} \}$.
\item $\mathcal{B}_3 = \{C \cup \{a,b,c \}: C \in \bigcup_{i=1}^{3} \mathcal{C}_i\}$.
\item $\mathcal{B}_4 = \{D \cup \{a,b,d \}: D \in \bigcup_{i=1}^{3} \mathcal{D}_i\}$.
\item $\mathcal{B}_5 = \{E \cup  \{c,d \}: E \in \bigcup_{i=1}^{3} \mathcal{E}_i\}$.
\end{list}
\begin{spec_splicing}\label{spec_splicing}
$(X \cup \{b,c,d\}, \bigcup_{i=1}^{5} \mathcal{B}_i)$ is a $(v + 3,k,4,6)$ cover.
\end{spec_splicing}
\begin{proof} Let $M$ be a $6$-subset of $X \cup \{b,c,d\}$:
\begin{list}{}{}
\item If $M \cap \{a,b,c,d\} = \emptyset$ or $\{a\}$ then $M$ is covered by some block $B \in \mathcal{B}_1$.
\item If $M \cap \{a,b,c,d\} = \{b\}$ or $\{c\}$ or $\{d\}$, then $M$ is covered by some block $B \in \mathcal{B}_1 \cup \mathcal{B}_2$.
\item If $\vert M \cap \{a,b,c,d\} \vert = 2$, then there exists $P \in\{\{a,b,c\}, \{a,b,d\},\{c,d\}\}$ such that $P \supseteq (M \cap \{a,b,c,d\})$. Let $P = \{a,b,c\}$ (the cases when $P = \{a,b,d\}$ or $P = \{c,d\}$ are similar). From property \ref{pro1}, $M$ is covered by some block $B \in \mathcal{B}_3$.
\item If $M \cap \{a,b,c,d \} = \{a,b,c \}$, then $M$ is covered by some block $B \in \mathcal{B}_3$.
\item If $M \cap \{a,b,c,d \} = \{a,b,d \}$, then $M$ is covered by some block $B \in \mathcal{B}_4$.
\item If $M \cap \{a,b,c,d \} = \{a,c,d \}$, then $M$ is a set of the form $\{a,c,d,x,y,z \}$ where $\{x,y,z \}$ is any triple of $ X \setminus \{a\}$. From property \ref{pro4}, it follows that there exists a block $T \in \bigcup_{i=1}^{3} \mathcal{C}_i \cup \mathcal{D}_i \cup \mathcal{E}_i$ such that $\vert T \cap \{x,y,z\} \vert \geq 2$. Moreover, since $\{a,c,d\}$ pairwise intersects in two points with $\{a,b,c\}, \{a,b,d \}$ and $\{c,d\}$, it follows that $ \vert M \cap B \vert \geq 4$  for some $B \in \bigcup_{i=3}^{5} \mathcal{B}_i$, $B \supset T$. The same arguments apply to the case in which $M \cap \{a,b,c,d \} = \{b,c,d \}$.
\item If $M \cap \{a,b,c,d\} = \{a,b,c,d\}$ then $M$ is clearly covered by some block $B \in \bigcup_{i=3}^{4} \mathcal{B}_i$.
\end{list}
The pair $(X \cup \{b,c,d \}, \bigcup_{i=1}^{5} \mathcal{B}_i)$ is therefore a $(v + 3,k,4,6)$ cover. \end{proof}

This construction implies the following upper bound formula:
\begin{equation*} %\label{splice6}
\begin{split}
C(v + 3,k,4,6) & \leq 3 \left\lfloor (k/v)
C(v,k,4,6)\right \rfloor + C(v,k,4,6) \\
& \quad + C(v_{1,1}+v_{2,1},k-3,2) + C(v_{1,2}+v_{2,2},k-3,2) \\
& \quad + C(v_{1,1}+v_{2,2},k-3,2) + C(v_{1,2}+v_{2,1},k-3,2)  \\
& \quad + C(v_{1,1}+v_{1,2},k-2,2) + C(v_{2,1}+v_{2,2},k-2,2) \\
& \quad + 2 C(v_{3,1}+v_{3,2},k-3,2) + C(v_{3,1}+v_{3,2},k-2,2),
\end{split}
\end{equation*}
where $\sum_{i=1}^{3} \sum_{j=1}^{2} {v}_{i,j} = v-1$ and $k \geq 5$.

\bigskip
By proceeding almost identically to the construction %for $(v + 3,k,4,6)$ covers 
presented above, we can build a $(v + 3,k,4,5)$ cover starting from a $(v,k,4,5)$ cover and obtain the following upper bound  for $C(v+3,k,4,5)$:
\begin{equation*}\label{splice5}
\begin{split}
C(v + 3,k,4,5) & \leq  3 \left\lfloor (k/v)
C(v,k,4,5)\right \rfloor + C(v,k,4,5) \\
& \quad + C(v_{1,1}+v_{2,1},k-3,2) + C(v_{1,2}+v_{2,2},k-3,2) \\
& \quad  + C(v_{1,1}+v_{2,2},k-3,2) + C(v_{1,2}+v_{2,1},k-3,2)  \\
& \quad + C(v_{1,1}+v_{1,2},k-2,2) + C(v_{2,1}+v_{2,2},k-2,2), \\
\end{split}
\end{equation*}
where $\sum_{i=1}^{2} \sum_{j=1}^{2} {v}_{i,j} = v-1$ and $k \geq 5$.

\bigskip
The point splicing constructions presented in this section and Section \ref{n-point} allow us to build new general covering designs based on ``smaller'' general covering designs. The nice feature of the point splicing construction of Section \ref{n-point} is that it can be used to obtain $t$-$(v,k,m, \lambda)$ general covering designs for arbitrary values of $v,k,t,m$ and $\lambda$.  

Starting from a $(v,k,4,m)$ cover with $k =6$ and $5 \leq m \leq 6$, two choices are available for the construction of a $(v+3,k,4,m)$ cover: the point splicing construction of Section \ref{n-point}, or the specific point splicing construction for $(v,k,4,m)$ covers presented in this section. Usually, if $k > 6$, the construction of Section \ref{n-point} performs better, but for $k = 6$ it is the point splicing construction of this section that gives the better results.
\section {Trapping-triples Construction for $(v,6,3,m)$ Covers} \label{constr63m}
In his paper \cite{decaen}, de Caen presented a construction for $(v,3,m)$ Tur\'an systems. It was based on the partition of a set $X$ into $m-1$ quasi-equal parts, that is, parts whose sizes pairwise differ by one unit at most. For $m =4$, de Caen's construction coincides with the one given by Tur\'an  in \cite{tur} who conjectured that it always produces optimal Tur\'an systems with $T(v,3,4)$ blocks. The conjecture has been shown to be true for $v \leq 13$ (\cite{spe}). Etzion et al. \cite{etz1} extended de Caen's construction to $(v,4,3,m)$ covers. It can be further extended to $(v,6,3,m)$ covers as follows:
Let $X$ be a set of $v$ elements and $X_0, X_1, \ldots, X_{m-2}$ be a partition of $X$ into $m-1$ quasi-equal parts. For $i = 0,1, \ldots, m-2$, let $(X_i,\mathcal{B}_i)$ be a
$(v_i,2,1)$ covering design with $w_i$ blocks $B_i^1, B_i^2, \ldots, B_i^{w_i}$, where $w_i = \lceil v_i/2 \rceil$, and let us select $h_i$
$(v_i,4,2)$ covering designs $(X_i,\mathcal{A}^1_i), (X_i, \mathcal{A}^2_i), \ldots, (X_i,\mathcal{A}^{h_i}_i)$  such that $(X_i, \bigcup_{j=1}^{h_i} \mathcal{A}^j_i)$ is a $(v_i,4,3)$ covering design, where $h_i = w_{(i + 1)\bmod(m-1)}$. Let us define
%, denoted by $ \mathcal{Q}_i$,
%forms the collection of
%blocks of a 
\[
\mathcal{C}_i = \bigcup_{j=1}^{h_i} \mathcal{A}^j_i \{B^j_{(i + 1)\bmod(m-1)}\},
\]
where $i =0,1,\ldots, m-2$.
\begin{teo63m}\label{teo63m}
$(X, \bigcup_{i=0}^{m-2}\mathcal{C}_i )$ is a $(v,6,3,m)$ cover.
\end{teo63m}
\begin{proof} Let us analyze how a given $m$-subset $M$ of $X_0 \cup X_1 \cup \ldots \cup X_{m-2}$ is covered in three points by some block $C$ in $\mathcal{C}_0 \cup \mathcal{C}_1 \cup \ldots \cup \mathcal{C}_{m-2}$:
\begin{enumerate}
\item Let $\vert M \cap X_l \vert \geq 3$ for some $l \in \{0,1, \ldots, m-2 \}$. From the definition of $\bigcup_{j=1}^{h_l} \mathcal{A}^j_l$,  it follows that for some $p \in \{1,2,\ldots,h_l\}$ there exists a block $A \in \mathcal{A}^p_l$ such that $\vert M \cap A \vert \geq 3$. This implies $\vert M \cap C \vert \geq 3$ where $ C = A \cup B^p_{(l + 1) \bmod(m-1)}$.
\item Let $\vert M \cap X_i \vert \leq 2$, for $0 \leq i \leq m-2$. Then there exists at least one and at most $\left \lfloor \frac{m}{2} \right \rfloor$ different parts $X_{j_n}$ of $X$ such that $\vert M \cap \ X_{j_n} \vert = 2$. This implies that there must exist $l \in \{0, 1, \ldots, m-2\}$ such that $\vert M \cap X_l \vert = 2$ and $\vert M \cap X_{(l + 1) \bmod(m-1)}\vert \geq 1$. From the definition of $h_l$ and $\mathcal{B}_{(l + 1) \bmod (m-1)}$, it follows that there exists $p \in \{1,2,\ldots,h_l\}$ such that $\vert M \cap B^p_{(l + 1) \bmod(m-1)} \vert \geq 1$. Since $(X_l, \mathcal{A}^p_l)$ is a $(v_l,4,2)$ covering design, there must exist a block $A \in \mathcal{A}^p_l$ such that $\vert M \cap \ A \vert = 2$ and therefore $\vert M \cap C \vert \geq 3$ where $ C = A \cup B^p_{(l + 1) \bmod(m-1)}$.
\end{enumerate}
We have therefore shown that, for any $m$-subset $M$ of $X$, there exists a block $C \in \bigcup_{i=0}^{m-2}\mathcal{C}_i$  such that $\vert M \cap C \vert \geq 3$. That is, $(X, \bigcup_{i=0}^{m-2}\mathcal{C}_i )$ is a $(v,6,3,m)$ cover. 
\end{proof}
From Theorem \ref{baker1976} and Theorem \ref{teirlinck1994} we can derive a general upper bound formula for $(v,6,3,m)$ covers.
\begin{genteo63m}\label{genteo63m}
Let $n$ be any positive integer. For $v = 4^n$ or $v = 2p^n + 2$ with $p \in \{7,31,127\}$, the following inequality holds:
\begin{displaymath}
C((m-1)v,6,3,m)\leq \frac{v^2(v-1)(m - 1)}{24}.
\end{displaymath}
\end{genteo63m}
\begin{proof}
Let $X$ be a set of $(m-1)v$ points. Let $v = 4^n$ or $v = 2p^n + 2$ where $n$ is a positive integer and $p \in \{7,31,127 \}$. For $i= 0,1,\ldots, m-2$: \\
Let $X_i$ be a part of $X$ and $ \vert X_i \vert= v$. Let $B^1_i, B^2_i, \ldots, B^w_i$ be a partition of $X_i$ where $\vert B^1_i \vert = \vert B^2_i \vert = \ldots = \vert B^w_i \vert = 2$. Therefore $w = \frac{v}{2}$.

From Theorem \ref{baker1976} and Theorem \ref{teirlinck1994} it follows that there exists a Steiner quadruple system $(X_i, \mathcal{A}_i)$ which is 2-partitionable. This implies that the collection $\mathcal{A}_i$  of blocks can be partitioned into $r$ parts $\mathcal{A}_i^1, \mathcal{A}_i^2, \ldots, \mathcal{A}_i^r$, each of which is the collection of blocks of a Steiner system $S(v,4,2)$. The value of $r$ is
\begin{displaymath}
\frac{v(v-1)(v-2)}{4 \cdot 3 \cdot 2} \cdot \frac{4 \cdot 3}{v(v-1)} = \frac{v-2}{2},
\end{displaymath}
where $\frac{v(v-1)(v-2)}{4 \cdot 3 \cdot 2}$ is the number of blocks of an SQS($v$) and $\frac{v(v-1)}{4 \cdot 3}$ is the number of blocks of an $S(v,4,2)$.
Let $\mathcal{A}_i^h$ be an additional collection of blocks such that $\mathcal{A}_i^h = \mathcal{A}_i^r$ and $h= r+1$. Clearly, $h = \frac{v}{2} = w$.

We now have all the elements to apply the construction presented in Theorem \ref{teo63m} which develops as follows:
\begin{displaymath}
\begin{split}
\left\vert \bigcup_{i=0}^{m-2}\mathcal{C}_i \right\vert & = \left\vert \bigcup_{j=1}^{h} \mathcal{A}^j_0 \{B^j_1\} \cup \bigcup_{j=1}^{h} \mathcal{A}^j_1 \{B^j_2\} \cup \ldots\cup \bigcup_{j=1}^{h} \mathcal{A}^j_{m-2} \{B^j_0\} \right\vert \\
& = (m - 1)\left(\frac{v(v-1)(v-2)}{4 \cdot 3 \cdot 2} + \frac{v(v-1)}{4 \cdot 3} \right) \\
& = (m - 1)\left(\frac{v(v-1)(v-2) + 2v(v-1)}{24}\right)\\
& = (m - 1) \left(\frac{v(v-1)(v-2+2)}{24}\right)\\
& = \frac{v^2(v-1)(m - 1)}{24}. \qedhere
\end{split}
\end{displaymath}
\end{proof}
As a consequence of Theorem \ref{genteo63m}, the following upper bound on the minimum size of $(3v,6,3,4)$ covers can be stated:
\begin{genteo634}\label{genteo634}
Let $n$ be any positive integer. For $v = 4^n$ or $v = 2p^n + 2$ with $p \in \{7,31,127\}$,
\begin{displaymath}
C(3v,6,3,4)\leq \frac{v^2(v-1)}{8}.
\end{displaymath}
\end{genteo634}
\section {Trapping-quadruples Constructions}\label{constrk4_m}
In the following section we present a construction for $(v,k,4,6)$ covers and sufficient conditions for its application will be discussed. Then, by requiring additional conditions to be satisfied, a construction for $(v,k,4,5)$ covers will be derived.
\subsection{Construction of $(v,k,4,6)$ Covers}\label{seck46}
Let $X$ be a set of $v$ elements, $v$ even, and $X_1$, $X_2$ be a partition of
$X$ into two equal parts. Let $n =\frac{v}{2}$. Moreover, let $k$ be an even number, $k \geq 4$ and $h = \frac{k}{2}$. Suppose there exists a resolvable $(n,h,2)$
covering design with $p$ parallel classes, $p \leq 5$. Let $\mathcal{P}_1, \mathcal{P}_2, \ldots, \mathcal{P}_p$ be the parallel classes defined on $X_1$ and $\mathcal{R}_1, \mathcal{R}_2, \ldots, \mathcal{R}_p$ be the parallel classes defined on $X_2$. For $i =1,2$, let $(X_i,\mathcal{B}_i)$ be an $(n,k,4)$ covering design. We assume therefore that $n \geq k > h$. Under this assumption, Theorem \ref{haemer_gen} implies $p\geq 3$. Define
\begin{displaymath}
\mathcal{B}= \mathcal{B}_1 \cup \mathcal{B}_2 \cup \bigcup_{i=1}^{p} \mathcal{P}_i\mathcal{R}_i.
\end{displaymath}
\begin{teok46}\label{teok46}
$(X, \mathcal{B})$ is a $(v,k,4,6)$ cover.
\end{teok46}
\begin{proof}  Let us analyze how a $6$-subset $M$ of $X_1 \cup X_2$ is covered in 4 points by some block $B \in \mathcal{B}$:
\begin{enumerate}
\item \label{case_4_covering_up}$M \in [6,0] \cup [5,1] \cup [4,2]$. Then there exists a block $B \in \mathcal{B}_1$ such that $\vert B \cap M \vert \geq 4$ since $(X_1, \mathcal{B}_1)$ is an $(n,k,4)$ covering design.
\item \label{case_4_covering_down}$M \in [0,6] \cup [1,5] \cup [2,4] $. Then there exists a block $B \in \mathcal{B}_2$ such that $\vert B \cap M \vert \geq 4$ since $(X_2, \mathcal{B}_2)$ is an $(n,k,4)$ covering design too.
\item $M \in [3,3]$. Let $T = M \cap X_1$ and $S = M \cap X_2$.
\begin{enumerate}
\item \label{case_contain}Suppose that the triple $T$ is contained in a block $P$ of a parallel class $\mathcal{P}_i$, $i \in \{1,2,\ldots, p \}$. Then from the definition of $\mathcal{R}_i$ it follows that there exists a block $R \in \mathcal{R}_i$ such that $\vert R \cap M \vert \geq 1$. This implies $\vert M \cap B \vert  \geq 4$ where $B = P \cup R$. We can proceed symmetrically when the triple $S$ is contained in a block of a parallel class $\mathcal{R}_j$, $j \in \{1,2,\ldots, p \}$.
\item \label{case_not_contain}Suppose instead that $T$ is not contained in any block of any class $\mathcal{P}_i$, $1 \leq i \leq p$, and $S$ is not contained in any block of any class $\mathcal{R}_j$, $1 \leq j \leq p$. Then for some $i_1,i_2, i_3 \in \{1,2,\ldots, p\}$, where $i_1 < i_2 < i_3$, and for some $j_1, j_2, j_3 \in \{1,2,\ldots, p\}$, where $j_1 < j_2 < j_3$, there must exist $I_1 \in \mathcal{P}_{i_1}$, $I_2 \in \mathcal{P}_{i_2}$, $I_3 \in \mathcal{P}_{i_3}$, $J_1 \in \mathcal{R}_{j_1}$, $J_2 \in \mathcal{R}_{j_2}$ and $J_3 \in \mathcal{R}_{j_3}$ such that $\vert T \cap I_l \vert = \vert S \cap J_l \vert = 2$, for $1 \leq l \leq 3$. This is because the pairs in $T$ (and the pairs in $S$) pairwise intersect in one point and cannot be contained in different blocks of a same parallel class by definition. Since $p \leq 5$, there must exist $y,z \in \{1,2,3 \}$ such that $i_y = j_z$. This implies $\vert (I_y \cup J_z) \cap M \vert = 4$.
\end{enumerate}
\end{enumerate}
We have shown that, for any $6$-subset $M$ of $X$, there exists a block $B \in \mathcal{B}$ such that $\vert M \cap B \vert \geq 4$. Hence $(X, \mathcal{B})$ is a $(v,k,4,6)$ cover. \end{proof}

Under the conditions of the construction presented in this section, we have
\begin{displaymath}
C(v,k,4,6) \leq 2C(v/2,k,4) + \frac{pv^2}{k^2}.
\end{displaymath}
\begin{remarkupper646}\label{remarkupper646}
From the construction mentioned above, we deduce that it is not always true that a given $6$-subset $M$ of $X_1 \cup X_2$, $M \in [4,2] \cup [2,4]$, is covered in four points by some block $B \in \bigcup_{i=1}^{p} \mathcal{P}_i\mathcal{R}_i$, but it is true if the size of each parallel class is less than four. Let us investigate the reason. Let $M \in [4,2]$ (the case when $M \in [2,4]$ can be dealt with in a similar way) and suppose that the size of each parallel class is $q < 4$. For $i = 1,2, \ldots, p$, the four points of the quadruple $M \cap X_1$ cannot lie in four different blocks of $\mathcal{P}_i$ (as the size of each class is less than four) and therefore $\vert (M \cap X_1) \cap P \vert \geq 2$ for some block $P \in \mathcal{P}_i$. On the other hand, there exists a parallel class $\mathcal{R}_j$, for some $j \in \{1,2, \ldots, p \}$, which contains a block $R$ such that $ \vert M \cap R \vert = 2$. This implies that, for some $P \in \mathcal{P}_j$, we have $\vert M \cap (P \cup R) \vert \geq 4$ and the above-mentioned construction can be improved by replacing $\mathcal{B}_1$ and $\mathcal{B}_2$ with the collections $\mathcal{C}_1$ and $\mathcal{C}_2$ of two $(n,k,4,5)$ covers $(X_1,\mathcal{C}_1)$ and $(X_2,\mathcal{C}_2)$. This improvement implies the following better upper bound for $(v,k,4,6)$ covers:
\begin{displaymath}
C(v,k,4,6) \leq 2C(v/2,k,4,5) + \frac{pv^2}{k^2}.
\end{displaymath}
\end{remarkupper646}
\subsection{\textbf {Construction of $(v,k,4,5)$ Covers}}
Let us consider again the construction presented in Section \ref{seck46} but instead of requiring that the number of parallel classes be $p \leq 5$, we require that the size of the parallel classes be $q = 2$.
\begin{teok45}\label{teok45}
$(X, \mathcal{B})$ is a $(v,k,4,5)$ cover.
\end{teok45}
\begin{proof} Let us analyze how a $5$-subset $M$ of $X_1 \cup X_2$ is covered in 4 points by some block $B \in \mathcal{B}$:
\begin{enumerate}
\item $M \in [5,0] \cup [4,1] \cup [1,4] \cup [0,5]$. Then $M$ is covered by some block $B \in  \mathcal{B}_1 \cup \mathcal{B}_2$ for the same considerations made in Theorem \ref{teok46}, points \ref{case_4_covering_up} and \ref{case_4_covering_down}.
\item $M \in [3,2] \cup [2,3]$. Let $\vert M \cap X_1 \vert = 3$ (the case when $\vert M \cap X_2 \vert = 3$ can be dealt with in a similar way). For $i = 1,2, \ldots, p$, the three points of the triple $M \cap X_1$ cannot lie in three different blocks of $\mathcal{P}_i$ (as the size of each class is less than three). By similarly following the same arguments made in Remark \ref{remarkupper646}, we deduce that for some $j \in \{1,2, \ldots, p \}$ there exist $P \in \mathcal{P}_j$ and $R \in \mathcal{R}_j$ such that
$\vert M \cap (P \cup R) \vert \geq 4$. \qedhere 
\end{enumerate}
\end{proof}
Under the conditions of the construction yielding Theorem \ref{teok45}, we have
\begin{displaymath}
C(v,k,4,5) \leq 2C(v/2,k,4) + 4p.
%\frac{pv^2}{k^2}
\end{displaymath}
Now, let us note that $C(k,k,4)$ is trivially equal to $1$ and that $C(3k,2k,4,5)=3$ \cite{burg}. Moreover, since $C(vm,km,t) \leq C(v,k,t)$ \cite{gordon1}, we  have $C(16m,8m,4) \leq C(16,8,4) = 30$ \cite{gordon}. These facts, combined with Theorem \ref{haemer}, Theorem \ref{teok46}, Remark \ref{remarkupper646} and Theorem \ref{teok45}, lead to the following upper bounds for covers:

\begin{teoupperformulak46}\label{teoupperformulak46} For $k \geq 2$, we have:
\begin{enumerate}
\item $C(4k, 2k, 4 ,5) \leq 14$ if $k$ is even; \label{case 3_2}
\item $C(4k, 2k, 4 ,5) \leq 18$ if $k$ is odd; \label{case 4_2}
\item $C(6k, 2k, 4 ,6) \leq 42$ if $k \equiv 0 \pmod 3$; \label{case 4_3}
\item $C(6k, 2k, 4 ,6) \leq 51$ if $k \equiv 1,2 \pmod 3$; \label{case 5_3}
\item $C(8k, 2k, 4 ,6) \leq 140$ if $k \equiv 0 \pmod 4$. \label{case 4_4}
\end{enumerate}
\end{teoupperformulak46}
Points \ref{case 3_2} and \ref{case 4_2} of Theorem \ref{teoupperformulak46} derive from Theorem \ref{haemer} and Theorem \ref{teok45}. Points \ref{case 4_3} and \ref{case 5_3} from Theorem \ref{haemer}, Theorem \ref{teok46} and Remark \ref{remarkupper646}. Point \ref{case 4_4} from Theorem \ref{haemer} and Theorem \ref{teok46}.

Here below some examples follow, where $p$ indicates the number of parallel classes and $q$ the size of each of them:
\begin{list}{}{}
\item $k = 3$, $p = 4$ and $q = 3$. In this case, the resolvable $(9,3,2)$ covering design with 12 blocks from Theorem \ref{haemer} is the well-known resolvable Steiner system $S(9,3,2)$. From Theorem \ref{teoupperformulak46} point \ref{case 4_3}, we have $C(18,6,4,6) \leq 42$, which matches the current best known upper bound for $C(18, 6, 4 ,6)$\footnote{Bertolo et al. (cf. \cite{bert}) applied an analogous technique to the case when $v=19$, $t = 4$ and $k=m=6$ from which the bound $C(18,6,4,6) \leq 42$ can be derived almost straightforwardly. However, they did not investigate the conditions under which the technique generalizes to other values of $v,k$ and $m$, nor did they determine upper bounds of the kind presented in Theorem \ref{teoupperformulak46}.}\cite{li}.
\item $k = 2$, $p = 5$ and $q = 3$. In this case, from Theorem \ref{haemer}, we have the $1$-factors of the $1$-factorization of the complete graph $K_6$. From Theorem \ref{teoupperformulak46} point \ref{case 5_3}, we have $C(12,4,4,6) \leq 51$, which matches the best upper bound for $C(12,4,4,6)$ \cite{gordon}, (i.e. for the Tur\'an number $T(12,4,6)$ and therefore for the covering number $C(12,8,6)$ as well).
\item $k = 2$, $p = 3$ and $q = 2$. In this case, from Theorem \ref{haemer}, we have the $1$-factors of the $1$-factorization of the complete graph $K_4$. From Theorem \ref{teoupperformulak46} point \ref{case 3_2}, we obtain $C(8,4,4,5) \leq 14$. Since an SQS$(8)$ exists, we have $C(8,4,3) = T(8,4,5) = C(8,4,4,5) = 14$. It is worth noting that in this case, from the construction yielding Theorem \ref{teok45}, blocks and coblocks not only have the same size but are identical: the constructed  $(8,4,5)$ Tur\'an system and its complement, a Steiner system $S(8,4,3)$, are the same design.
\end{list}
\subsection{A Construction of $(v,5,4,6)$ Covers}\label{sec5_4_6}
Let $X$ be a set of $v$ elements. Partition $X$ into 4 quasi-equal parts $X_0$, $X_1$, $X_2$ and $X_3$. For $i = 0,1,2,3$, let $v_i = \vert X_i\vert $. For $0 \leq i \leq 3$ and $1\leq j\leq h$, where $h$ is some positive integer, let $(X_i, \mathcal{A}_{i,j})$ be a $(v_i,3,2)$ covering design and $(X_i, \mathcal{B}_{i,j})$ be a $(v_i,2,1)$ covering design. For $0 \leq i \leq 3$, assume that $\mathcal{A}_{i,1} \cup \mathcal{A}_{i,2} \cup \ldots \cup \mathcal{A}_{i,h} = \binom{X_i}{3}$ and that $\mathcal{B}_{i,1} \cup \mathcal{B}_{i,2} \cup \ldots \cup \mathcal{B}_{i,h} = \binom{X_i}{2}$.
% is a $(v_i,3,2)$ covering design $(X_i, \mathcal{B}_{i})$.
For $i=0,1$, let $(X_i, \mathcal{D}_i)$ be a $(v_i,5,4,5)$ cover. Finally, for $i=2,3$, let $(X_i, \mathcal{E}_i)$ be a $(v_i,5,4)$ covering design.
Define
\[
\mathcal{C} = \bigcup_{i=0}^1\mathcal{D}_i \cup \mathcal{E}_{i+2} \cup \bigcup_{i=0}^3 \bigcup_{j=1}^h \mathcal{A}_{i,j}\mathcal{B}_{{(i+1) \bmod 4},j} \cup \bigcup_{i=0}^1\bigcup_{j=1}^h \mathcal{A}_{i,j}\mathcal{B}_{{(i+2)},j}.
\]
\begin{teo546}\label{teo546}
$(X, \mathcal{C})$ is a $(v,5,4,6)$ cover.
\end{teo546}
\begin{proof} Let us analyze how a given 6-subset $M$ of $X$ is covered in 4 points by some block $C \in \mathcal{C}$:
\begin{enumerate}
\item $\vert M \cap X_r \vert \geq 5$, for some $r \in \{0,1\}$. Then there exists a block $C \in \mathcal{D}_r$ such that $\vert C \cap M \vert \geq 4$ as $(X_r, \mathcal{D}_r)$ is a $(v_r,5,4,5)$ cover.
\item $\vert M \cap X_r \vert \geq 4$, for some $r \in \{2,3\}$. Then there exists a block $C \in \mathcal{E}_r$ such that $\vert C \cap M \vert \geq 4$ as $(X_r, \mathcal{E}_r)$ is a $(v_r,5,4)$ covering design.
\item $\vert M \cap X_r \vert \geq 2$, $\vert M \cap X_s \vert \geq 2  $ for some $r,s \in \{0,1,2,3 \}$, $r < s$. Observe that $\mathcal{A}_{z,1} \cup \mathcal{A}_{z,2} \cup \ldots \cup \mathcal{A}_{z,h} = \binom{X_z}{3}$ and $\mathcal{B}_{z,1} \cup \mathcal{B}_{z,2} \cup \ldots \cup \mathcal{B}_{z,h} = \binom{X_z}{2}$ for $z \in \{r,s\}$. Then, if $r = 0$ and $s=3$, it follows that for some $j \in \{1,2, \ldots, h\}$ there exists a block $C \in \mathcal{A}_{s,j}\mathcal{B}_{r,j}$ such that $\vert C \cap M \vert \geq 4$; otherwise, it follows that for some $j \in \{1,2, \ldots, h\}$ there exists a block $C \in \mathcal{A}_{r,j}\mathcal{B}_{s,j}$ such that $\vert C \cap M \vert \geq 4$.
\item $\vert M \cap X_r \vert = 3$ for some $r \in \{0,1,2,3 \}$, and $\vert M \cap X_s \vert = 1$ for all $s \in \{0,1,2,3 \}$, $s \neq r$. Since $\mathcal{A}_{r,1} \cup \mathcal{A}_{r,2} \cup \ldots \cup \mathcal{A}_{r,h} = \binom{X_r}{3}$, it follows that for some $j \in \{1, 2, \ldots, h \}$ there exists a block $A \in \mathcal{A}_{r,j}$ such that $\vert A \cap M \vert = 3$. On the other hand, there exists a block $B \in \mathcal{B}_{{(r+1) \bmod 4},j}$ such that $\vert B \cap M \vert = 1$. This implies $\vert (A \cup B) \cap M \vert = 4$.
%, where $C = A \cup B$, $A \cup B \in  \mathcal{A}_{r_j}\mathcal{B}_{{(r+1) \bmod 4},j}$.
\item $\vert M \cap X_r \vert =4$ for some $r \in \{0,1\}$, and $\vert M \cap X_{s_1} \vert =1$, $\vert M \cap X_{s_2} \vert =1$ for some $s_1, s_2 \in \{0,1,2,3 \}$, $s_1 < s_2$,  $r \neq s_1, r \neq s_2$. If $ r= 0$, it follows that $s_1 \in \{1,2\}$. Then, for some $j \in \{1,2, \ldots, h\}$, there exists a block $C \in \mathcal{A}_{r,j}\mathcal{B}_{{s_1},j}$ such that $\vert C \cap M \vert = 4$. Otherwise, if $r= 1$, it follows that $s_2 \in \{2,3\}$ and again, for some $j \in \{1,2, \ldots, h\}$, there exists a block $C \in \mathcal{A}_{r,j}\mathcal{B}_{{s_2},j}$ such that $\vert C \cap M \vert = 4$.
\end{enumerate}
The block design $(X, \mathcal{C})$ is indeed a $(v, 5, 4, 6)$ cover.
\end{proof}
Fort and Hedlund \cite{fort} proved that an optimal $(v,3,2)$ covering design has $\lceil \frac{v}{3} \lceil \frac{v-1}{2} \rceil \rceil$ blocks. Moreover,  for $v \geq 8$, we have $\mu(v,3) = v-2$ (\cite{etz}). With regard to $(v,2,1)$ covering designs, we have $\lambda(v,2) = \mu(v,2) = v-1$ if $v$ is even, otherwise, if $v$ is odd, $\lambda(v,2) = v-2$ and $\lambda(v,2) \lceil \frac{v}{2} \rceil = \frac{v(v-1)}{2} - 1$ (\cite{etz}). These results, combined with  Theorem \ref{teo546}, lead to the following upper bound for $(4v,5,4,6)$ covers:
\[
C(4v,5,4,6) \leq 2(C(v,5,4,5)+C(v,5,4)) + 3v(v-1) \Big\lceil \frac{v}{3} \Big\lceil \frac{v-1}{2} \Big\rceil \Big\rceil.
\]
\section{Conclusions}
Improving upper bounds on the minimum size of general covering designs is a challenging problem. In order to obtain good upper bounds, combinatorial constructions involving unions and intersections of different kind of combinatorial designs are very effective and we think it is worth keeping exploiting these techniques.

If a $2$-partitionable SQS$(v)$ exists, then, as Theorem \ref{genteo63m} shows, $C((m-1)v,6,3,m) \leq 24^{-1}(v^3 - v^2) (m-1)$. Apart from this implication on general covers, it would be very interesting to see in the future new results on the existence of $2$-partitionable Steiner quadruple systems.

\bigskip\noindent \textbf{Acknowledgement}
 \newline The author would like to thank Prof. Peter J. Cameron for his helpful comments.

\end{document}